# Comment: Classifier Technology and the Illusion of Progress

**Robert A. Stine**

It is my pleasure to contribute to the discussion of this paper. David Hand has the credibility one needs to write such an article and not have it dismissed out of Hand. Along with publishing numerous papers and books on classification and data mining, he "works in the trenches" with real data. His contributions to credit modeling are particularly well known and respected, and his knowledge of that domain reaches far deeper into the substance than the casual illustration often chosen to show off a new methodology. He is a fascinating lecturer and I have learned a great deal by listening carefully to his ideas. When he writes that claims of the superiority of neural networks and support vector machines "fail to take account of important aspects of real problems," I have to stop and think about my own research and experiences.

The thrust of Hand's paper is the argument that most recent developments in classification, say anything since Fisher's linear discriminant function, offer little benefit in practice. The mismatch between theory and practice dwarfs incremental claims for superiority established in theorems. For instance, theory that shows that a support vector machine classifies better than a simple linear model is an "illusion," bordering on sophistry.

I have a great deal of sympathy for this point of view, but I doubt that many statisticians will change what they do after reading this paper. I agree with many of his criticisms, but I am already in the choir. I suspect that it will take quite a bit more to convince others, particularly along the lines of proposals for what ought to be done. Consider the impact of Tukey's "The future of data analysis" (Tukey, 1962). After chastising the field for its preoccupation with "optimization in terms of a precise, similarly inadequate criterion," Tukey proposed alternatives, including exploratory data analysis and robust methods. Forty years later, Hand's criticisms echo his concerns.

*Robert A. Stine is Professor, Department of Statistics, The Wharton School, University of Pennsylvania, Philadelphia, Pennsylvania 19104-6340, USA (e-mail: stine@wharton.upenn.edu).*

Hand presents a range of criticisms of modern classifiers. I find it useful to organize my discussion by grouping them into two clusters:

- Creeping incrementalism
- Square pegs in round holes.

Let me start with the first of these.

*Creeping incrementalism.* Hand argues that concerns for optimality emphasize tiny improvements that are dwarfed by other issues in real applications. He argues that the first predictor or the most simple of models finds most of the structure. Adding bells and whistles contributes little more than complex window dressing, and the advantages are illusions that disappear during the application. The argument is analogous to saying that linear Taylor series make pretty good approximations to most functions; generally, you do not need those messy, higher order terms. I certainly agree that simple models— or at least simple methodologies—take you a long way. Dean Foster and I wrote a paper to make just this point when mining financial data: with a few adjustments, stepwise linear regression can predict bankruptcy as well as elaborate trees (Foster and Stine, 2004).

A convincing argument for preferring simpler models requires careful discussions of applications. Given the depth of his experience, I had expected Hand to offer a rich portfolio of examples that demonstrate the failures of complex models. Instead, he relies more on an idealized example (one of equally correlated predictors) and a summary of fitted models to selected data sets from the repository at UC Irvine. One has to be careful basing arguments on made-up examples, because it is too easy to turn the examples around. With equally correlated predictors, the first one or two predictors capture most of the signal, with diminishing benefits left to the others. Although I have had similar experiences modeling real data, it is all too easy to make up normal models in which later variables appear to explain the most variation. For example, define

$$
\begin{aligned}
X_1 &= \tau Y + \varepsilon_1 + \varepsilon_2, \\
X_2 &= \tau Y + \varepsilon_1 - \varepsilon_2, \\
X_3 &= \tau Y - \varepsilon_1 + \varepsilon_3, \\
X_4 &= \tau Y - \varepsilon_1 - \varepsilon_3,
\end{aligned}
\quad \text{where } Y, \varepsilon_i \overset{\text{i.i.d.}}{\sim} N(0,1). \tag{1}
$$

Each predictor $X_j$ has equal correlation $\tau$ with $Y$ and the predictors have a block structure. In this setting, what happens as we greedily expand the regression model is shown in Table 1. With $\tau = 0.25$, we have superadditive growth in the fit of the model: the addition of a subsequent predictor adds more to the model than any predecessor. I am not claiming that this example is more natural than the one in the paper. That is not the point. The point is that, separated from a real application, it is easy to construct examples that support any argument. What matters is what is useful in practice, and we need to see more evidence from real applications to appreciate the flaws of complex models.

I think that one needs to "go easy" when it comes to criticism of the use of statistical inference to judge improvements in a model. Inferential statistics concerns the separation of even a little signal from noise. This perspective is ideally suited to applications in traditional science. Discovery of statistically significant anomalies from the standard theory *is* important. A statistically significant anomaly, even a small one, cannot be dismissed as random variation and leads to revisions of the current theory. However, there needs to be a current theory in the first place. Without an established point of reference, the yardstick used to gauge improvements should be different. Most real applications lack such a benchmark and resemble an entirely new domain. When I was first learning about the connection between statistics and information theory, I was interested in the use of statistical models for data compression. (Think of tools used to compress the files on your computer disk.) Early on, improvements to algorithms for data compression regularly brought reductions of 20 or 30% in the amount of disk space required to store a data file. As the area matured, the gains got smaller and issues of statistical significance became relevant. Statistical significance in this context amounts to resolving whether you can save two or three more bits!

It is also important to establish what it means for a model to be better than another. Statistical significance offers one scale, but it may be poorly suited to the task. Finding an acceptable alternative can be particularly hard (e.g., in the social sciences), but is often easy in business. In business, improvements generally get measured in dollars, and statistical significance seldom guarantees much in the way of economic benefits. This point needs to be stressed as prominently and concretely as possible. Hand discusses the choice of the loss function used to judge classifiers and rightfully criticizes the casual use of error rates. Unfortunately, the survey of fitted models summarized in his Table 1, however, compares error rates. Who is to say that a small improvement in predictive accuracy is not valuable? Consider the data set "Segmentation" in the first row of his Table 1. Perhaps the reduction in the error rate from 0.083 to 0.014 is worth quite a lot of money. Without deeper insights into these applications, I cannot judge whether the improvements are impressive or unimportant. I doubt that enough is known about these applications to set costs, but perhaps Hand could offer other examples from his own experience in which the costs are known.

*Square peg in the round hole.* Statistics has rightly been criticized for often devoting too much energy to unrealistic problems. As Tukey pointed out, "Better to have an approximate answer to the right question than the exact answer to the wrong question." Knowing the right question, however, often means knowing more about the application than most of us get from clients. In working with banks on credit modeling, the proprietary nature of their business makes it nearly impossible for them to be able to disclose enough for me to think that I am answering the right question. That does not mean that I have stopped trying, but it gets painful to jam your foot in the door over and over. It can be a lot more satisfying to prove a theorem or write code for a new algorithm.

TABLE 1

| Number of predictors | Explained variation | | |
|---|---|---|---|
| | | $\tau = 0.25$ | $\tau = 0.5$ |
| 0 | | 0 | 0 |
| 1 | $\frac{1}{1+2/\tau^2} =$ | 0.03 | 0.11 |
| 2 | $\frac{1}{1+1/(2\tau^2)} =$ | 0.11 | 0.33 |
| 3 | $\frac{1}{1+1/(5\tau^2)} =$ | 0.24 | 0.55 |
| 4 | | 1 | 1 |

Another reason for solving the wrong problem is that by the time one has the data and builds a model, the problem has changed. I would push to the front of the line to agree with Hand that changes in the underlying population pose a serious problem. This problem is particularly acute in business because of its competitive environment. If a company builds a model that produces a change in its behavior (such as a better way to evaluate the risk of loans that it makes), you can be sure that the competition will react and change as well.

I recently had a first-hand experience with this type of problem. The task was to help a company improve the methods that it uses to evaluate prospective employees. Based on attributes known at the time of an application, we developed a classifier that was able to identify those most likely to succeed. The usual sorts of validation exercises showed that the effects we found were real, at least for the population represented by our data. As pointed out by Hand in Section 3, it takes a long time to get the data needed for this type of modeling. In our case, we had to wait and see which employees succeeded before we got the response. The delay was two years. By the time that the company tried to use the model, the economy had changed and the nature of the people applying for jobs had shifted. In fact, because we identified certain factors as important, the company changed the way that it collected these factors, rearranging the application form to emphasize the presentation of the key questions. I have little doubt that the revised questions measure different things than those used to build our model. Our model was a disappointment, but then I doubt that any model would have handled these disruptions.

I owe a favorite example of how the use of a model changes the population to Professor Hand. Suppose we are building a model to score the credit-worthiness of our customers. We discover that customers who, like me, drive white cars are poor risks. As a result, we stop offering loans to those driving white cars. Now think about what happens in several years when it is time to refresh the scoring model. By this time, none of our customers drives a white car, so this characteristic no longer appears to be a risk factor. Our successor will have to learn this all over again—that is, if these drivers have not changed their color preference. In the utopian world of repeated sampling from the population, these things do not happen. The population does not change because you start to use a model.

*What next?* Einstein once remarked, "Everything should be made as simple as possible, but not simpler." Given a preference for simple models, I would very much like for Hand to offer some guidance *suited to applications* on how one is supposed to decide whether it is useful to look for more structure. If not by the ruler given by statistical inference, then how? In my toy example, the sum $X_1 + X_2 + X_3 + X_4$ predicts $Y$ perfectly. What should we do, however, when we have a wide data set with relatively few cases and 1000 predictors? How would we know to try the sum of them all as a predictor? Stepwise methods that build up models are good at finding subadditive models, but superadditive structures are difficult to identify. Similarly, we have methods that capture nonlinear features in data, but how are we to know whether to try them? If we only look for simple models, then we will always find simple models. To find nonlinearities requires that we entertain models that allow them. For example, our regression model for predicting bankruptcy uses interactions that, in effect, segment the population. Without them, the predictions were much less able to predict bankruptcies and left a lot of money on the table (Foster and Stine, 2004).

Professor Hand has had more experience with the challenges of dealing with real applications than most statisticians. I would be very interested in his approach to deciding when additions to a simple model *are* worthwhile. Similarly, what are his thoughts on methods to assess population drift? Certainly, statisticians have been concerned about population drift for a long time. For example, consider the article by Brown, Durbin and Evans (1975) on detecting changes in a linear model, Kalman filters that explicitly model an evolving state variable or models for evolutionary time series dating back to Priestley (1965). Do these fail in practice?